\documentclass[11pt,twocolumn,twoside]{article}
\usepackage{fully3d}
\usepackage{amsmath}
\usepackage{amsfonts}
\usepackage{amssymb}
\usepackage{xcolor}
\usepackage{siunitx}
\usepackage{caption}
\usepackage{subcaption}
\usepackage{cleveref}
\usepackage{algorithm2e}

\let\Vec\mathbf
\begin{document}

\title{Software Implementation of the Krylov Methods Based Reconstruction for the 3D Cone Beam CT Operator}

\author[1]{Vojtěch Kulvait}
\author[1]{Georg Rose}

\affil[1]{Institute for Medical Engineering and Research Campus STIMULATE, University of Magdeburg, Magdeburg, Germany}

\maketitle
\thispagestyle{fancy}

\begin{customabstract}
Krylov subspace methods are considered a standard tool to solve large systems of linear algebraic equations in many scientific disciplines such as image restoration or solving partial differential equations in mechanics of continuum. In the context of computer tomography however, the mostly used algebraic reconstruction techniques are based on classical iterative schemes. In this work we present software package that implements fully 3D cone beam projection operator and uses Krylov subspace methods, namely CGLS and LSQR to solve related tomographic reconstruction problems. It also implements basic preconditioning strategies. On the example of the cone beam CT reconstruction of 3D Shepp-Logan phantom we show that the speed of convergence of the CGLS clearly outperforms PSIRT algorithm. Therefore Krylov subspace methods present an interesting option for the reconstruction of large 3D cone beam CT problems.
\end{customabstract}

\section{Introduction}
Krylov subspace methods, see \cite{JorgLiesen2012}, are attractive in the context of the solving tomographic problems as they do not require direct storage of the system matrix of the corresponding problem. This is especially true for 3D tomographic problems, where the size of the system matrix prohibits its efficient storage and the computation of the projection and back-projection operator on GPUs or dedicated hardware might be orders of magnitude faster than using stored precomputed matrix.

While most of the literature regarding algebraic reconstruction is based on Kaczmarz algorithm and related classical iterative schemes, see \cite{Gregor2008}, there is also a growing body of the works on application of various Krylov subspace methods. Important topics are optimal preconditioning strategies, see \cite{Cools2015} and the enforcing of properties such as non-negativity of the solution \cite{Gazzola2017}.

Until circa 2010 the direct application of the Krylov subspace methods to the 3D cone-beam CT operator (CBCT), was very rare and the works regarding Krylov methods were solving just smaller 2D problems. The block-wise algorithm to divide the tomographic reconstruction into the smaller sub-problems to apply CGLS and LSQR was proposed in \cite{Qiu2012}. Currently the implementation of CGLS for CBCT operator can be found in the two MATLAB based frameworks, see \cite{Aarle2015, Biguri2016}. The Split Bregman algorithm for CBCT TV norm minimization using Krylov BiCGStab was published in \cite{Molina2018,Vorst1992}. There is however still lack of works to systematically study Krylov subspace methods or to compare their performance with the classical ART based approaches.

In this work we present the software package, which contains an open-source C++ and OpenCL implementation of the Krylov subspace methods for the CBCT reconstruction. We show that these methods poses a viable option for a fast and accurate reconstruction on recent GPU hardware. 

Moreover on simple test problem using CBCT we compare performance of PSIRT, an advanced technique based on ART like algorithms, to the CGLS, implementation of conjugate gradients on the normal equations. We show that convergence speed of the Krylov method is much faster and to achieve the same accuracy, we need a lot less iterations. \bigskip
\section{Materials and Methods}

Cone beam CT operator can be understood as a sparse matrix $\Vec{A} \in \mathbb{R}^{m \times n}$ acting on discretized attenuation data $\Vec{x} \in \mathbb{R}^{n}$ in the volume of interest to produce projection data $\Vec{b}\in\mathbb{R}^{m}$. In the matrix form
\begin{equation}
    \Vec{b} = \Vec{A} \Vec{x}.
\end{equation}
As the matrix $\Vec{A}$ is non square and often over-determined with $m>n$, we typically solve the least-squares problem by means of normal equations to find attenuation $\Vec{x}$ such that
\begin{equation}
   \Vec{A}^\top \Vec{A} \Vec{x} = \Vec{A}^\top \Vec{b}.
\end{equation}

The matrix $\Vec{A}^\top \Vec{A}$ is symmetric, positive definite and it is possible to apply directly method of conjugate gradients on such system. Direct method to do so is referred as CGLS. There is also LSQR, mathematically equivalent method with improved numerical stability, see \cite{Paige1982,Reichel2008}. For the sake of completeness, we include here the iterative scheme of CGLS implemented in our software as Algorithm~\ref{algorithm}. We restructured the algorithm in a way that at the end of each iteration we compute the update of $\Vec{x}$ and postpone the update of the residuals to the beginning of the next iteration. By doing so, we save one projection and backprojection at the end of the algorithm.
\begin{algorithm}[t]
\SetKwInOut{Input}{input}

 \Input{Projection data $\Vec{b}$, initial vector $\Vec{x}_0$, relative discrepancy tolerance $\mathrm{ERR}$, maximum number of iterations $K$.}
 \Begin{
    allocate $\Vec{x}$, $\Vec{d}_x$ and $\Vec{r}_x$\;
    allocate $\Vec{e}_b$ and $\Vec{p}_b$\;
    $\mathrm{NB}_0 = \|\Vec{b}\|_2$\;
    $\Vec{x} = \Vec{x}_0$\;
    $\Vec{p}_b = \Vec{A} \Vec{x}$\;
    $\Vec{e}_b = \Vec{b} - \Vec{p}_b$\;
    $\Vec{r}_x = \Vec{A}^\top \Vec{e}_b$\;
    $\Vec{d}_x = \Vec{r}_x$\;
    $\mathrm{NR2}_\mathrm{old} = \|\Vec{r}_x\|_2^2$\;
    $\Vec{p}_b = \Vec{A} \Vec{d}_x$\;
    $\mathrm{NP2} = \|\Vec{p}_b\|_2^2$\;
    $\alpha =  \mathrm{NR2}_\mathrm{old} / \mathrm{NP2}$ \;
    $\Vec{x} = \Vec{x} + \alpha \Vec{d}_x$\;
    $\Vec{e}_b = \Vec{e}_b - \alpha \Vec{p}_b$\;
    $\mathrm{NB} = \|\Vec{e}_b\|_2$\;
    $i=1$\;
    \While{$\mathrm{NB}/\mathrm{NB}_0 > \mathrm{ERR} \And i < \mathrm{K}$ }{
    $\Vec{r}_x = \Vec{A}^\top \Vec{e}_b$\;
    $\mathrm{NR2}_\mathrm{now} = \|\Vec{r}_x\|_2^2$\;
    $\beta = \mathrm{NR2}_\mathrm{now} / \mathrm{NR2}_\mathrm{old}$\;
    $\Vec{d}_x = \Vec{d}_x + \beta \Vec{r}_x$\;
    $\mathrm{NR2}_\mathrm{old} = \mathrm{NR2}_\mathrm{now}$\;
    $\Vec{p}_b = \Vec{A} \Vec{d}_x$\;
        $\mathrm{NP2} = \|\Vec{p}_b\|_2^2$\;
    $\alpha =  \mathrm{NR2}_\mathrm{old} / \mathrm{NP2}$ \;
    $\Vec{x} = \Vec{x} + \alpha \Vec{d}_x$\;
    $\Vec{e}_b = \Vec{e}_b - \alpha \Vec{p}_b$\;
    $\mathrm{NB} = \|\Vec{e}_b\|_2$\;
     $i = i+1$\;
 }
     }
\KwResult{Vector $\Vec{x}$, number of iterations $i$, final norm of discrepancy $\mathrm{NB}$.}
 \caption{CGLS with delayed residual computation.}
 \label{algorithm}
\end{algorithm} \section{Software}

The software package was developped in C++ and OpenCL.  The project implements various CBCT projection and back-projection operators, namely Siddon ray-caster \cite{Siddon1985}, footprint methods \cite{Long2010} and also so called Cutting voxel projector. Cutting voxel projector uses the volume of the cuts of the voxels by the rays to particular pixel for the computation of projections, details are yet to be published. From the reconstruction techniques, the software contains CGLS and LSQR implementation with the option of basic Jacobi preconditioning and Tikhonov regularization. It is possible to select initial vector or guess of the solution, e.g. the result of analytical reconstruction or apriori knowledge. It is also possible to do a off-center reconstruction, where the volume to reconstruct can be positioned outside the center of rotation. For the purpose of the comparison of the different CBCT reconstruction methods, two classical ART like methods, SIRT and PSIRT \cite{Gregor2008} are also implemented.

The package also contains methods to project volumes or backproject projections without reconstruction. This could be useful when e.g. simulating acquisition of the given volume with particular geometry setting of a concrete CT device.

The program is distributed under the terms of GNU GPL3 license and its Git repository is available at \url{https://github.com/kulvait/kct_cbct}. The results presented were obtained using git commit f2bf01a.

\section{Results}
Here we present a test to compare convergence of the Krylov method (CGLS) with the classical scheme (PSIRT) when applied on CBCT operator. We have chosen 3D Shepp-Logan phantom with the 256x256x52 voxels of the dimensions $\SI{0.86}{\milli \meter} \times \SI{0.86}{\milli \meter} \times \SI{3.44}{\milli \meter}$. 

The geometry configuration is similar to the clinical C-Arm CT systems for the brain tomographic scanning, where the distance from the source to the isocenter is $\SI{749}{\milli \meter}$ and the distance from source to the detector is $\SI{1198}{\milli \meter}$. Detector matrix have pixel size of $\SI{0.616}{\milli \meter} \times \SI{0.616}{\milli \meter}$. The trajectory consists of $496$ scanning angles.  

We have first projected the 3D Shepp-Logan phantom using this geometry to obtain projection data. We have used for the projections and the reconstructions the implementation of TT projector and backprojector, see \cite{Long2010}. We run the tests to compare classical method PSIRT to Krylov subspace method CGLS, both implemented in our software. The tests were performed on computer with the AMD Ryzen 7 1800X and GPU Vega 20 Radeon VII with 16GB HBM2 Memory and 1TB/s bandwidth. Projectors and backprojectors are implemented in OpenCL and the computations were run on the GPU. Speed of the both methods in terms of average time per iteration was comparable, circa $\SI{12.9}{\second}$. To compare the speed of convergence, we have measured relative norm of discrepancy of the solution
\begin{equation}
e=\frac{\|\Vec{A}\Vec{x} - \Vec{b}\|_2}{\|\Vec{b}\|_2}. \label{eq:discrepancy}
\end{equation}
during the iterative process. 

\begin{figure*}[t]
  \centering
  \includegraphics[width=1.0\textwidth]{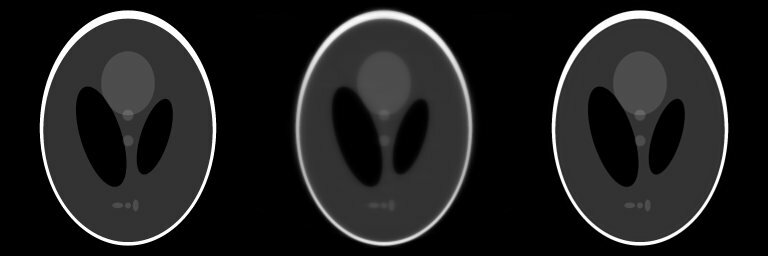}
  \caption{Center slice of the test volume of 256x256x52 voxels, all images have the same window [0,1], ground truth data on the left. PSIRT reconstruction after 40 iterations in the middle. CGLS reconstruction after 40 iterations on the right.}
  \label{psirtcgls}
\end{figure*}

Initially we run fixed number of $40$ iterations of the both methods, the relative norm of discrepancy \eqref{eq:discrepancy} after $40$ iterations was then $e_{\mathrm{CGLS}}= 0.18\%$ versus $e_{\mathrm{PSIRT}}= 3.64\%$. The visualization in Figure~\ref{psirtcgls} shows that PSIRT reconstruction is still blurry while CGLS has converged without visible problems.

\begin{figure}[t]
  \centering
  \includegraphics[width=1.0\columnwidth]{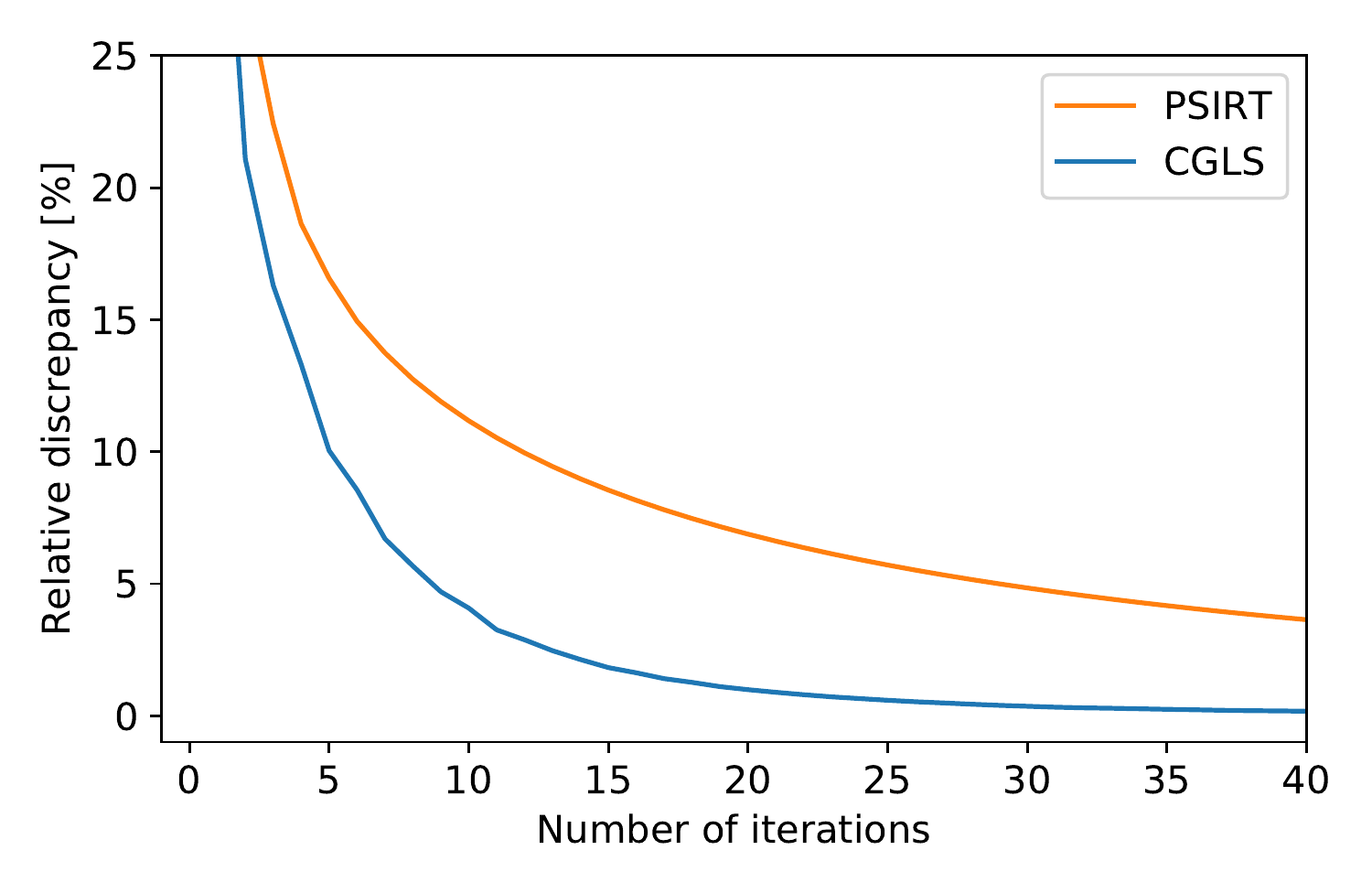}
  \caption{Comparison of the speed of convergence in terms of the relative norm of discrepancy \eqref{eq:discrepancy} for CGLS and PSIRT for the test CBCT problem, 40 iterations.}
  \label{convergence_cut}
\end{figure}

\begin{figure}[t]
  \centering
  \includegraphics[width=1.0\columnwidth]{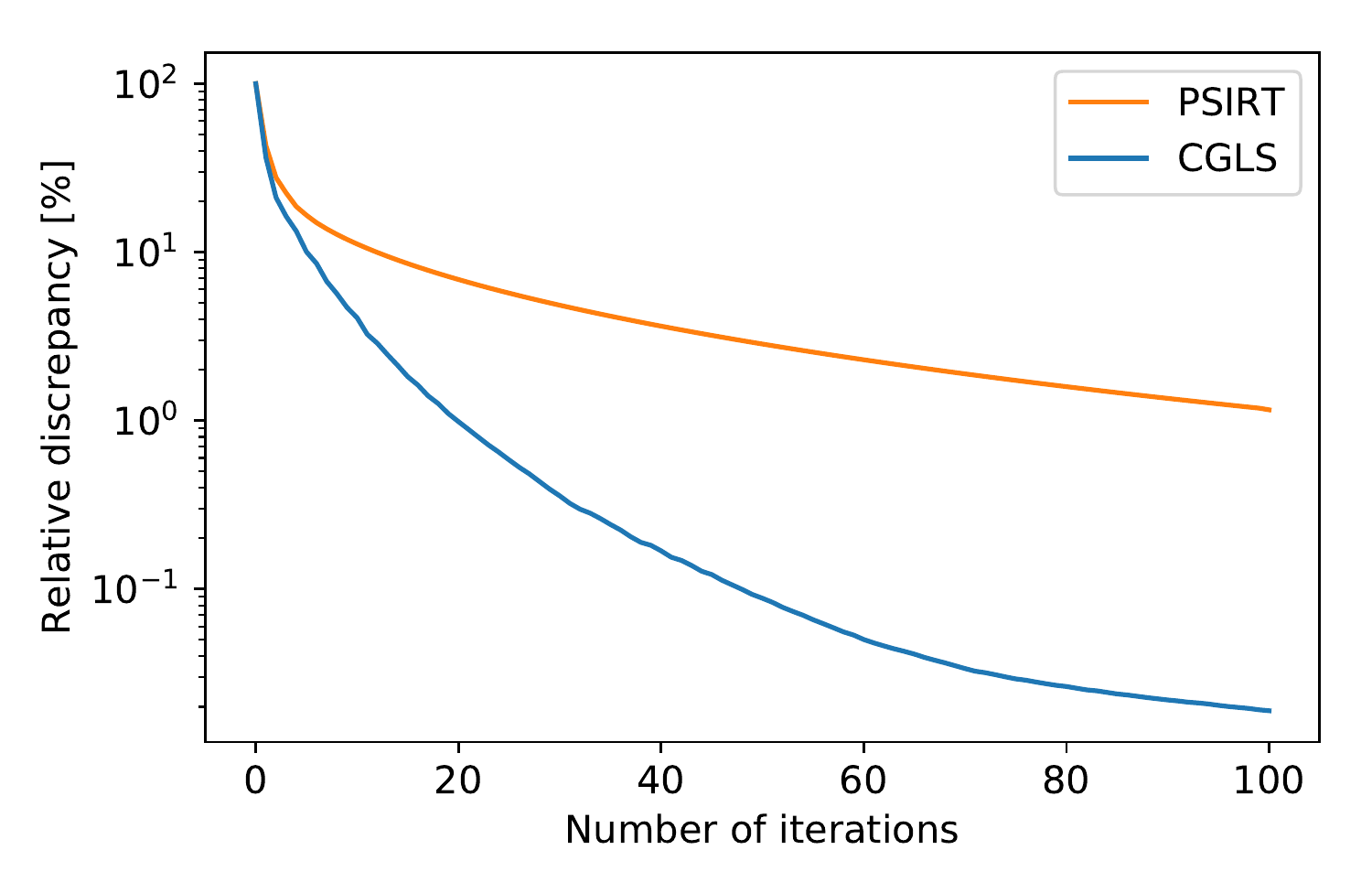}
  \caption{Comparison of the speed of convergence in terms of the relative norm of discrepancy \eqref{eq:discrepancy} for CGLS and PSIRT on the test CBCT problem. Y axis in logarithmic scale,  100 iterations.}
  \label{convergence}
\end{figure}

Second, we test how many iterations every method needs to achieve norm of discrepancy under $e < 1\%$. For PSIRT it is $N=112$ iterations compared to CGLS with $N=20$. This means that the CGLS is circa five times faster than PSIRT when we would like to achieve the same accuracy. In \cref{convergence_cut,convergence} can be found the graphs comparing the speed of convergence for both methods that illustrates clear advantage of the CGLS for the test problem.

 \section{Discussion}
The library we present in this article was developed to test different approaches to implement cone beam CT operator and to test various reconstruction techniques. Focus was on solving problems that arise from the model based perfusion reconstructions, where we can have negative values in reconstructed volumes and classical solvers might often diverge. The accuracy is often prioritized over speed and many double precision computations could be substituted by float precision counterparts to increase speed. From the nature of the open-source solution and the C++ modular implementation it is straightforward to extend the program for the particular reconstruction algorithm or preconditioning strategy.

On a test problem we have demonstrated that the CGLS has much faster convergence characteristic than PSIRT. It can be seen from \cref{psirtcgls,convergence,convergence_cut}. We have tested also SIRT as alternative classical reconstruction method, but its convergence was slower when compared to PSIRT.

When using Krylov methods, it is not well addressed how to enforce conditions such as non-negativity of the solution. While in the context of classical ART, this is an easy task, for Krylov methods we can not simply apply the box conditions as the regularized solution would not fit to the underlying Krylov subspace and we would lose the convergence properties of the Krylov methods. To address this multiple schemes for restarted methods within this framework were introduced, see \cite{Gazzola2017}. Implementing these techniques is also one of further goals for the development of presented software package. Main obstacle to do so is a larger memory footprint of such methods as compared to CGLS we usually need to store more vectors of underlying Krylov space.

The method PSIRT could be used with the box conditions as for the test problem values lies within a range $[0,1]$. The relative norm of discrepancy after 40 iterations was $e_{\mathrm{PSIRT}}= 3.64\%$ without box conditions and $e_{\mathrm{PSIRT}}= 3.62\%$ with them. Therefore for simplicity, we report only the results of PSIRT without box conditions as they did not have any meaningful effect on the convergence. In real world applications the upper bound can be hard to estimate.

The memory footprint of SIRT and CGLS is the same. Our implementation of CGLS, see Algorithm \ref{algorithm}, needs to store three times the volume data and three times the right hand side data. It is necessary for storing the residual and Krylov subspace vectors on which we project the error. Our implementation of SIRT has the same memory footprint as we need to store the vectors of the row and column sums of the system matrix and update vectors. In PSIRT we do not store column sums vector so we need to store right hand side data only twice. However, current GPU hardware provides enough memory such that the memory footprint of these methods is not an issue in practical applications even for much bigger problems than the one presented.

There is also implemented Jacobi preconditioning strategy in the software. It is most simple preconditioning, where we approximate the matrix of the normal equations by its diagonal. It seems however that this strategy alone does not work very well especially due to the presence of very small diagonal values in system matrix on the cone boundary. Therefore better preconditioning strategies have to be found in order to speed the convergence.

Although LSQR should in theory be numerically more stable than CGLS, from our experiments for the tomographic reconstruction the two methods are producing practically identical results. So the type of instability that makes LSQR numerically more stable method is probably not present in a typical CT data. Due to smaller memory footprint, using CGLS might therefore be preferred.

In the CGLS algorithm the discrepancy and residual vectors are not computed directly but they are iteratively updated. Therefore it is proposed to reorthogonalize this vector once in $k$ iterations to avoid accumulation of errors. When we applied such scheme, we figured out that the difference between iteratively updated discrepancy and discrepancy computed from solution vector is less than $0.0001\%$ for $10$ iterations. Therefore we omit this reorthogonalization step in a default configuration.

 \section{Conclusion}
Analytical reconstruction methods are still a gold standard in a CT reconstruction. Main argument for using them is their speed. As Krylov subspace methods provide much faster convergence when compared to the ART like methods, together with the hardware speed improvements, their wider application could change the speed narrative to widely adopt algebraic reconstruction techniques in practice.

The results show on a phantom problem very high advantage of the CGLS over PSIRT in terms of convergence and they shall be validated for practical problems and other setups. We believe that there will be still very strong advantage of Krylov solvers in practice. Further development should focus on adopting a good preconditioning strategies for Krylov solvers as it has potential to further reduce the run time.

\section{Remarks}
Just before the conference, we managed to achieve a significant speedup of some projectors and backprojectors implemented in the \url{https://bitbucket.org/kulvait/kct_cbct}, especially the so-called Cutting voxel projector, which will be introduced in a separate article. Because of this speedup, significantly better projection, backprojection and reconstruction times were presented on the poster. The main subject of this contribution, the convergence speed of Krylov methods as a function of the number of projections and backprojections, is not affected. However, it may significantly improve the usability of our software for fast algebraic reconstruction of moderate sized problems in minutes. Last stable commit as of writing this sentence is \texttt{0d7d6}.

After presenting a poster at the Fully3D 2021 conference, we received very important feedback. Namely, Simon Rit, the lead developer of the Reconstruction Toolkit, see \cite{Rit2014}, mentioned that methods derived from  the Kaczmarz algorithm using the entire projection and backprojection operator at once are slow compared to ordered subsets methods, see \cite{Kong2009}. This also applies to the PSIRT algorithm, which does not use ordered subsets and which we compare Krylov methods with. We take this criticism very seriously and have started working on implementing OS algorithms in our package to be able to reliably compare the methods. Unfortunately, this comparison is not yet complete and cannot be presented here. It is also worth mentioning that, compared to ordered subset schemes, it is easy to add L2 regularization directly in the problem formulation for Krylov methods. As recently shown, L2 regularization can be a faster alternative to L1 regularization, such as TV norm minimization, with similar results, see \cite{Bilgic2014}.

The previous question however, leads to the following consideration: if ordered subsets methods have yielded significant speedups of classical CT reconstruction schemes by using only a subset of the rows of the CT operator matrix in each step, can a similar approach be used for Krylov methods? A naive implementation of methods like CGLS with this strategy would very likely suffer from a rapid loss of orthogonality and convergence speed. On the other hand, it is possible that orthogonalization with respect to a larger number of vectors proportional to the number of subsets could stabilize the method. Although this would lead to longer recurrences and a larger memory footprint of the algorithm, it seems promising to investigate such methods in the future.  
\section*{Acknowledgments}
 \textit{The work of this paper is partly funded by the Federal Ministry of Education and Research within the Research Campus STIMULATE under grant number 13GW0473A.} 

\printbibliography

\end{document}